\documentclass[12pt]{article}
\usepackage{mathrsfs}
\usepackage{CJK}
\usepackage{bbm}
\usepackage{stmaryrd}
\usepackage{shortvrb,psfrag}
\usepackage{enumerate}
\usepackage{amsmath}
\usepackage{cases}
\usepackage{amsfonts}
\usepackage{latexsym,amssymb,amsmath,mathrsfs,amsbsy, amsthm}
\usepackage[dvips]{graphicx}
\usepackage[usenames]{color}
\usepackage{xspace,colortbl}
\usepackage{pdfsync}
\usepackage[dvipdfm,
pdfstartview=FitH, CJKbookmarks=true, bookmarksnumbered=true,
bookmarksopen=true,
citecolor=blue
,colorlinks=true
,pdfborder=1
,pdfstartpage=2
,pdfstartview=Fit]{hyperref}
\allowdisplaybreaks

\newtheorem{thm}{Theorem}[section]
\newtheorem{lem}[thm]{Lemma}

\newtheorem{rem}{Remark}
\newtheorem{cor}[thm]{Corollary}

\newtheorem{prop}[thm]{Proposition}

\newcommand{\beq}{\begin{equation}}
\newcommand{\eeq}{\end{equation}}
\newcommand{\ben}{\begin{eqnarray}}
\newcommand{\een}{\end{eqnarray}}
\newcommand{\beno}{\begin{eqnarray*}}
\newcommand{\eeno}{\end{eqnarray*}}

\begin{document}
\title{{\bf Liouville-type theorems for the planer stationary MHD equations with growth at infinity
}
\\Wendong WANG
\\[2mm]
{\small $ ^\ddag$ School of Mathematical Sciences, Dalian University of Technology, Dalian 116024,  China}\\
{\small E-mail: wendong@dlut.edu.cn}\\[2mm]}

\date{\today}
\maketitle

% A critical Hardy inequality and its application to

\begin{abstract}
For the two dimensional stationary MHD equations, we prove that Liouville type theorems hold if the velocity is growing at infinity, where the magnetic field is assumed to be bounded under a smallness condition. The key point is to overcome the nonlinear terms, since no maximum principle holds for the MHD case with respect to the Navier-Stokes equations. As a corollary, we obtain that all the solutions of the 2D Navier-Stokes equations satisfying $\nabla u\in L^p(\mathbb{R}^2)$ with $1<p<\infty$ are constants, which is sharp since the same argument fails  in the case of $\nabla u\in L^\infty(\mathbb{R}^2)$.
\end{abstract}

{\small {\bf Keywords:} Liouville theorem, MHD equations, Navier-Stokes equations}

{\bf 2010 Mathematics Subject Classification:} 35Q30, 35B53, 76D03.

\setcounter{equation}{0}
\section{Introduction}

Consider the incompressible MHD equations on the whole space $\mathbb{R}^2$:
\begin{equation}\label{eq:MHD-2D}
\left\{\begin{array}{llll}
-\Delta u+u\cdot \nabla u+\nabla \pi=b\cdot\nabla b,\\
-\Delta b+u\cdot \nabla b=b\cdot\nabla u,\\
{\rm div }~ u=0,\quad {\rm div }~ b=0,
\end{array}\right.
\end{equation}
and the Dirichlet energy is defined as follows:
\ben\label{eq:energy bound-MHD-2D}
D(u,b)=\int_{R^2}|\nabla u|^2+|\nabla b|^2dx.
\een

%Here $\nu$ and $\mu$ denote viscosity and resistivity coefficients.

When $b=0$ in (\ref{eq:MHD-2D}), it follows that the 2D Navier-Stokes equations. Let us recall some known results on
 this issue. For example,
Gilbarg-Weinberger proved the above Liouville type theorem by assuming (\ref{eq:energy bound-MHD-2D}) in  \cite{GW1978}, where they made use of the fact that the vorticity function satisfies
a nice elliptic equation to which a maximum principle applies. The assumption on boundedness of the Dirichlet energy can be relaxed to $\nabla u\in L^p(\mathbb{R}^2)$ with some $p\in (\frac65,3]$, see Bildhauer-Fuchs-Zhang \cite{BFZ2013}.
If $u$ is bounded, a Liouville theorem being more in the spirit of the classical one for entire analytic functions
was obtained by Koch-Nadirashvili-Seregin-Sverak \cite{KNSS} as a byproduct of their
work on the non-stationary case. The above results also can be generalized to the shear thickening flows, for example see \cite{Fu2012exist,Fu2012Liou,FZ2012,ZG2013,JK2014,ZG2015}. The existence and asymptotic behavior of solutions in an exterior domain, for example see \cite{GNP1997,Ru2009,Ru2010,GG2011,PR2012,KPR2014,DI2017}.

Moreover, a velocity field $u$ satisfying the stationary Navier-Stokes equations on the
entire plane must be constant under the growth condition $\limsup |x|^{-\alpha}|u(x)|<\infty$ as
$|x|\rightarrow\infty$ for some $\alpha\in [0, 1/7)$, see Fuchs-Zhong \cite{FZ2011}. The component is improved to $\alpha<\frac13$, see Bildhauer-Fuchs-Zhang \cite{BFZ2013}.
More references, we refer to \cite{KNSS} and the references therein.

For the two dimensional stationary MHD equations, the similar Liouvile type theorems seem to be more difficult, since the maximum principle is
not available to the best of my knowledge. In \cite{WW}, the author and Y. Wang obtained some Liouvile type theorems by assuming (\ref{eq:energy bound-MHD-2D}) or $u\in L^\infty$, where the smallness conditions of the magnetic field are added. Here we go on this topic in this direction. Since all the exact solutions of (\ref{eq:MHD-2D}) with $b=0$ we know are polynomials, it seems that the smooth solutions below linear growth are trivial. A natural question:

{\bf What happens if the velocity is growing at infinity?}

Note that the vorticity equations are as follows.
Let $w=\partial_2u_1-\partial_1 u_2$ and $h=\partial_2b_1-\partial_1 b_2$, then
\begin{equation}\label{eq:MHDv-2D}
\left\{\begin{array}{llll}
-\Delta w+u\cdot \nabla w=b\cdot\nabla h,\\
-\Delta h+u\cdot \nabla h=b\cdot\nabla w+ H
\end{array}\right.
\end{equation}
where
\beno
H=2\partial_2b_2(\partial_2u_1+\partial_1 u_2)+2\partial_1u_1(\partial_2b_1+\partial_1 b_2).
\eeno

The main difficulty comes from the terms $b\cdot\nabla w$, $H$ etc., which is not vanishing for any energy integration. That's why we have to assume the smallness of some norm of  $b$.
However, if the velocity is largely growing as in \cite{BFZ2013}, i.e. there exist two constants
$\alpha>0$ and $c_0>0$ such that
\ben\label{eq:growth at infty-u}
|u(x)| \leq c_0\left(1+|x| \right)^{\alpha},\quad \forall~x\in \mathbb{R}^2,
\een
it's more complicated in this case. In fact, as the same arguments in \cite{BFZ2013}, the term
\beno
C(q)R^{2\alpha}\int_{\mathbb{R}^2}h^2~ |b|^2w^{2q-4}(\eta^{2\ell})dx
\eeno
seems to be out of control
(see (\ref{eq:energy estimate-w}) in the second subsection). To overcome it, we
introduce the decay condition of $b$:
 \ben\label{eq:growth at infty-b}
|b(x)| \leq c_0\left(1+|x| \right)^{\beta},\quad \forall~x\in \mathbb{R}^2,
\een
where $\beta<0$, and consider the local energy estimate in an annular domain.

Next we state our first result:
\begin{thm}\label{thm:generalization of KNSS}
Let $(u,b,\pi)$ be a smooth solution of the 2D MHD equations (\ref{eq:MHD-2D}) defined over the entire plane satisfying the growth estimates (\ref{eq:growth at infty-u}) with $\alpha<\frac13$ and (\ref{eq:growth at infty-b}) with $\beta<-\alpha$. Then $u$ and $\pi$ are constants and $b\equiv0$, if there exists one positive number
$\varepsilon_0=\varepsilon_0(\alpha,\beta,c_0)$ such that there holds
\beno
\|b\|_{L^{1}(\mathbb{R}^2)}+\||h|^{\frac13}\|_{L^{1}(\mathbb{R}^2)}\leq \varepsilon_0.
\eeno
\end{thm}

\begin{rem}
The above result generalized the Liouville type theorem in
\cite{KNSS,FZ2011,BFZ2013} to the MHD case.
\end{rem}

%L^{\frac{q_0^2}{q_0-1}}(\mathbb{R}^2)\cap
It follows from the above theorem that
\begin{thm}\label{thm:generalization of KNSS'}
Let $(u,b,\pi)$ be a smooth solution of the 2D MHD equations (\ref{eq:MHD-2D}) defined over the entire plane satisfying the growth estimates $\nabla u\in L^{q_0}(\mathbb{R}^2)$ for $1< q_0<\infty$, $\nabla b\in  L^{\infty}(\mathbb{R}^2)$ and
\beno
\|b\|_{L^1(\mathbb{R}^2)}+\||h|^{1/3}\|_{L^1(\mathbb{R}^2)}\leq \varepsilon,
\eeno
where $\varepsilon$ is sufficiently small depending on $q_0,$ $\|\nabla u\|_{L^{q_0}}$ and $\|\nabla b\|_{L^\infty}$. Then $u$ and $\pi$ are constants and $b\equiv0$.
\end{thm}

When $b$ vanishes, the 2D Navier-Stokes equations follows from (\ref{eq:MHD-2D}).
\begin{cor}\label{cor:generalization of KNSS'}
Let $(u,\pi)$ be a smooth solution of the 2D NS equations  defined over the entire plane satisfying the growth estimates $\nabla u\in L^q(\mathbb{R}^2)$ for some $1<q<\infty$.
 Then $u$ and $\pi$ are constants.
\end{cor}

\begin{rem}
The above result generalized the Liouville type theorem by  Gilbarg-Weinberger in \cite{GW1978} for $q=2$. Moreover,  this is the best estimate in a sense, since there are counter-examples for  $\nabla u\in L^\infty(\mathbb{R}^2)$ (for example, the Couette flow $(x_2,0)$).
\end{rem}

%
%Recall a result of Gilbarg-Weinberger in \cite{GW1978}.
%\begin{lem}\label{lem:GW}
%Let $f=(f_1,f_2)\in C^1$ in $r>r_0$ and have finite Dirichlet integral
%\beno
%\int_{r>r_0}|\nabla f|^2dxdy<\infty.
%\eeno
%Then we have
%\beno
%\lim_{r\rightarrow\infty} \frac{1}{\ln r}\int_0^{2\pi}|f(r,\theta)|^2d\theta=0,
%\eeno
%and
%there is a sequence $\{r_n\}$, $r_n\in (2^n, 2^{n+1})$, such that
%\beno
%\lim_{n\rightarrow \infty}\frac{|f(r_n,\theta)|^2}{\ln r_n}=0
%\eeno
%uniformly in $\theta.$
%\end{lem}

%If, furthermore, we assume $\nabla f \in L^p(\mathbb{R}^2)$ for some $2< p < \infty$,
%then the above decay property can be improved to be point-wise uniformly.
%More precisely, we have
We need the following lemma in the proof.
\begin{lem}[Theorem II.9.1 \cite{Galdi}]
\label{lem:Galdi}
Let $\Omega \subset \mathbb{R}^2$ be an exterior domain and let
\beno
\nabla f \in  L^p (\Omega),
\eeno
for some $2< p < \infty$.
Then
\beno
\lim_{|x| \rightarrow \infty} \frac{|f(x)|}{|x|^{\frac{p-2}{p}}} = 0,
\eeno
uniformly.
\end{lem}

Throughout this article, $C(\alpha_1,\cdots,\alpha_n)$ denotes a constant depending on $\alpha_1,\cdots,\alpha_n$, which may be different from line to line.

\section{Proof of Theorem \ref{thm:generalization of KNSS}}

In this section, we are aimed to prove Theorem \ref{thm:generalization of KNSS} by following the same route in \cite{BFZ2013}. Different from the arguments in \cite{BFZ2013}, we consider the local energy estimates in an annular domain and obtain the $L^q$ estimates of the vorticity.

First, we prove the following proposition.
\begin{prop}\label{thm:generalization of KNSS2}
Let $(u,b,\pi)$ be a smooth solution of the 2D MHD equations (\ref{eq:MHD-2D}) defined over the entire plane satisfying the growth
estimates (\ref{eq:growth at infty-u}) with $0<\alpha<\frac13$ and (\ref{eq:growth at infty-b}) with $\beta<-\alpha$.
Then
%\beno
%\nabla u, \nabla b\in L^q(\mathbb{R}^2),
%\eeno
%and
\beno
\|\nabla u\|_{L^{2q}(\mathbb{R}^2)}+\|\nabla b\|_{L^{2q}(\mathbb{R}^2)}\leq C(\alpha,\beta,q,c_0)<\infty
\eeno
holds for any $q>q_0$, where
\ben\label{eq:q0}
q_0=\max\{\frac{2}{1-3\alpha}, \frac{-1}{\alpha+\beta},-\frac{1}{2\beta}\}.
\een
\end{prop}

{\bf Proof of Proposition \ref{thm:generalization of KNSS2}.}
Let $\eta(x)\in C_0^\infty(B_R)$ be a cut-off function on an annular domain with $0\leq \eta\leq 1$ satisfying
\begin{align*} \eta(x)=\left\{
\begin{aligned}
&1,\quad x\in B_R\backslash B_{R/2},\\
&0, \quad x\in B_{2R}^c\cup B_{R/4}.
\end{aligned}
\right. \end{align*}
Write $w^{2q}=(w^2)^{q}$. Then for $q\geq 2,\ell\geq q$, we have
\beno\label{eq:w estimate-mhd}
\int_{\mathbb{R}^2}w^{2q}\eta^{2\ell}dx\nonumber&=& \int_{\mathbb{R}^2}(\partial_2 u_1-\partial_1 u_2)w^{2q-2}w\eta^{2\ell}dx\nonumber\\
&=&  \int_{\mathbb{R}^2}(u_2,-u_1)\cdot \nabla [w^{2q-2}w\eta^{2\ell}]dx\nonumber\\
&\leq &(2q-1)\int_{\mathbb{R}^2}|u||\nabla w| w^{2q-2}\eta^{2\ell}dx+2\ell\int_{\mathbb{R}^2}|u|| \nabla \eta| |w|^{2q-1}\eta^{2\ell-1}dx\nonumber\\
&\leq&\frac12 \int_{\mathbb{R}^2}w^{2q}\eta^{2\ell}dx+C(q)\int_{\mathbb{R}^2}|u|^2|\nabla w|^2 w^{2q-4}\eta^{2\ell}dx\nonumber\\
&&+2\ell\int_{\mathbb{R}^2}|u|| \nabla \eta| |w|^{2q-1}\eta^{2\ell-1}dx
\eeno
Similarly, we have
\beno\label{eq:w estimate-mhd2}
\int_{\mathbb{R}^2}h^{2q}\eta^{2\ell}dx
&\leq&C(q)\int_{\mathbb{R}^2}|b|^2|\nabla h|^2 h^{2q-4}\eta^{2\ell}dx+4\ell\int_{\mathbb{R}^2}|b|| \nabla \eta| |h|^{2q-1}\eta^{2\ell-1}dx
\eeno
%and
%\ben\label{eq:w estimate-mhd3}
%\int_{R^2}h^2w^{2q-3}\eta^{2\ell}dx&=& \int_{R^2}(\partial_2 b_1-\partial_1 b_2)hw^{2q-3}\eta^{2\ell}dx\nonumber\\
%&\leq&\delta \int_{R^2}w^{2q}\eta^{2\ell}dx+C(\delta,q)\int_{R^2}|u|^2|\nabla w|^2 w^{2q-4}\eta^{2\ell}dx\nonumber\\
%&&+2\ell\int_{R^2}|u|| \nabla \eta| w^{2q-1}\eta^{2\ell-1}dx
%\een
Due to the growth estimates (\ref{eq:growth at infty-u}) and (\ref{eq:growth at infty-b}),  we have
\ben\label{eq:w estimate-mhd4}
&&\int_{\mathbb{R}^2}w^{2q}\eta^{2\ell}+h^{2q}\eta^{2\ell}dx\nonumber\\
&\leq&C(q)R^{2\alpha}\int_{\mathbb{R}^2}|\nabla w|^2 w^{2q-4}\eta^{2\ell}dx+C(q)R^{2\beta}\int_{\mathbb{R}^2}|\nabla h|^2 h^{2q-4}\eta^{2\ell}dx\nonumber\\
&&+ C(\ell)R^{\alpha-1}\int_{\mathbb{R}^2} |w|^{2q-1}\eta^{2\ell-1}dx+C(\ell)R^{\beta-1}\int_{\mathbb{R}^2}|h|^{2q-1}\eta^{2\ell-1}dx
\een

On the other hand,
multiply $\eta^{2\ell}w^{2q-4}w$ and $\eta^{2\ell}h^{2q-4}h$ on both sides of (\ref{eq:MHDv-2D}), and we have
\ben\label{eq:energy estimate-w}
I&\doteq&(2q-3)\int_{\mathbb{R}^2}|\nabla w|^2 w^{2q-4}\eta^{2\ell}dx\nonumber\\
&\leq&\frac{1}{2q-2}\int_{\mathbb{R}^2} w^{2q-2}\triangle (\eta^{2\ell})dx+\frac{1}{2q-2}\int_{\mathbb{R}^2} w^{2q-2}u\cdot \nabla(\eta^{2\ell})dx\nonumber\\
&&+\int_{\mathbb{R}^2} b\cdot \nabla h |w|^{2q-4}w(\eta^{2\ell})dx \nonumber\\
%&\leq&\frac{1}{2q-2}\int_{\mathbb{R}^2} w^{2q-2}\triangle (\eta^{2\ell})dx+\frac{1}{2q-2}\int_{\mathbb{R}^2} w^{2q-2}u\cdot \nabla(\eta^{2\ell})dx\nonumber\\
%&&-(2q-3)\int_{\mathbb{R}^2}h~ b\cdot \nabla w w^{2q-4}(\eta^{2\ell})dx -\int_{\mathbb{R}^2}w^{2q-3}h~ b\cdot \nabla (\eta^{2\ell})dx\nonumber\\
&\leq&\frac{1}{2q-2}\int_{\mathbb{R}^2} w^{2q-2}\triangle (\eta^{2\ell})dx+\frac{1}{2q-2}\int_{\mathbb{R}^2} w^{2q-2}u\cdot \nabla(\eta^{2\ell})dx\nonumber\\
&&+\frac12 I+C(q)\int_{\mathbb{R}^2}h^2~ |b|^2w^{2q-4}(\eta^{2\ell})dx -\int_{\mathbb{R}^2}w^{2q-4}w h~ b\cdot \nabla (\eta^{2\ell})dx
\een
and similarly
\ben\label{eq:energy estimate-h}
II&\doteq&(2q-3)\int_{\mathbb{R}^2}|\nabla h|^2 h^{2q-4}\eta^{2\ell}dx\nonumber\\
%&=&\frac{1}{2q-2}\int_{\mathbb{R}^2} h^{2q-2}\triangle (\eta^{2\ell})dx+\frac{1}{2q-2}\int_{\mathbb{R}^2} h^{2q-2}u\cdot \nabla(\eta^{2\ell})dx\nonumber\\
%&&+\int_{\mathbb{R}^2} b\cdot \nabla w ~h^{2q-3}(\eta^{2\ell})dx+\int_{\mathbb{R}^2} H ~h^{2q-3}(\eta^{2\ell})dx\nonumber\\
&\leq&\frac{1}{2q-2}\int_{\mathbb{R}^2} h^{2q-2}\triangle (\eta^{2\ell})dx+\frac{1}{2q-2}\int_{\mathbb{R}^2} h^{2q-2}u\cdot \nabla(\eta^{2\ell})dx\nonumber\\
&&-\int_{\mathbb{R}^2}  w ~|h|^{2q-4}h b\cdot \nabla(\eta^{2\ell})dx+\frac12 II+C(q)\int_{\mathbb{R}^2}w^2~ |b|^2h^{2q-4}(\eta^{2\ell})dx \nonumber\\
&&+C\int_{\mathbb{R}^2} |\nabla u||\nabla b||h|^{2q-3}(\eta^{2\ell})dx
\een
Then it follows from (\ref{eq:w estimate-mhd4}), (\ref{eq:energy estimate-w}) and  (\ref{eq:energy estimate-h}) that
\ben\label{eq:w estimate-mhd-general}
&&\int_{\mathbb{R}^2}w^{2q}\eta^{2\ell}+h^{2q}\eta^{2\ell}dx\nonumber\\
&\leq&C(q,\ell)R^{\alpha-1}\left(R^{\alpha-1}\int_{\mathbb{R}^2} w^{2q-2}(\eta^{2\ell-2})dx+R^{2\alpha}\int_{\mathbb{R}^2} w^{2q-2}(\eta^{2\ell-1})dx+\int_{\mathbb{R}^2} |w|^{2q-1}(\eta^{2\ell-1})dx\right)\nonumber\\
&&+C(q,\ell)\left(R^{2\alpha+2\beta}\int_{\mathbb{R}^2} h^2w^{2q-4}(\eta^{2\ell})dx+R^{2\alpha-1+\beta}\int_{\mathbb{R}^2} |w|^{2q-3}|h|(\eta^{2\ell-1})dx\right)\nonumber\\
&&+C(q,\ell)R^{\beta-1}\left(R^{\beta-1}\int_{\mathbb{R}^2} h^{2q-2}(\eta^{2\ell-2})dx+R^{\alpha+\beta}\int_{\mathbb{R}^2} h^{2q-2}(\eta^{2\ell-1})dx+\int_{\mathbb{R}^2} |h|^{2q-1}(\eta^{2\ell-1})dx\right)\nonumber\\
&&+C(q)\left(R^{4\beta}\int_{\mathbb{R}^2} w^2h^{2q-4}(\eta^{2\ell})dx+R^{-1+3\beta}\int_{\mathbb{R}^2} |h|^{2q-3}|w|(\eta^{2\ell-1})dx\right)\nonumber\\
&&+C(q)R^{2\beta} \int_{\mathbb{R}^2} |\nabla u||\nabla b||h|^{2q-3}(\eta^{2\ell})dx  =I_1+\cdots+I_5
\een

{\bf Estimate of $I_5$.}
For a smooth vector-valued function $F\in C_0^2(\Omega)$, by applying the Calder\'{o}n-Zygmund theory we have
\ben\label{eq:CZ estimate}
\|\nabla F\|_{L^q(\Omega)}\leq C(n,q)\left(\|{\rm div}~ F\|_{L^q(\Omega)}+\|\nabla\times F\|_{L^q(\Omega)}\right),
\een
since the following identity holds,
\beno
\triangle F=\nabla({\rm div}~F)-\nabla\times\nabla\times F.
\eeno
Hence, by choosing $F=u \eta^{\frac{\ell}{q}}$ or $b \eta^{\frac{\ell}{q}}$  we get
\beno
\left(\int_{\mathbb{R}^2}|\nabla u|^{2q}\eta^{2\ell}dx\right)^{\frac{1}{2q}}&\leq& C(q,\ell) \left(\int_{\mathbb{R}^2}| u|^{2q}\eta^{2\ell-2q}|\nabla\eta|^{2q}dx\right)^{\frac{1}{2q}}+ C(q,\ell)\left(\int_{\mathbb{R}^2}|w|^{2q}\eta^{2\ell}dx\right)^{\frac{1}{2q}}\\
&\leq& C(q,\ell)R^{-1+\alpha+\frac1q}+ C(q,\ell)\left(\int_{\mathbb{R}^2}|w|^{2q}\eta^{2\ell}dx\right)^{\frac{1}{2q}}
\eeno
%and due to $b\in L^\infty$
%\beno
%\left(\int_{\mathbb{R}^2}|\nabla b|^{2q}\eta^{2\ell}dx\right)^{\frac{1}{2q}}
%&\leq& C(q,\ell)R^{-1+\frac1q}+ C(q,\ell)\left(\int_{\mathbb{R}^2}|h|^{2q}\eta^{2\ell}dx\right)^{\frac{1}{2q}}
%\eeno
and
\beno
&&\int_{\mathbb{R}^2} |\nabla u||\nabla b||h|^{2q-3}(\eta^{2\ell})dx\\
&\leq& C(q,\ell)R^{\frac1q}\left(R^{\alpha+\frac1q-1}  +\left(\int_{\mathbb{R}^2}|w|^{2q}\eta^{2\ell}dx\right)^{\frac{1}{2q}} \right)\cdot\left(R^{\beta+\frac1q-1}  +\left(\int_{\mathbb{R}^2}|h|^{2q}\eta^{2\ell}dx\right)^{\frac{1}{2q}} \right) \\
&&\cdot \left(\int_{\mathbb{R}^2}h^{2q}\eta^{2\ell}dx\right)^{\frac{2q-3}{2q}}.
%&\leq& C(\delta,q,\ell)R^{\frac{2q}{3}(-1+\alpha+\beta+\frac3q)}+ C(\delta,q,\ell)R^{q(-1+\alpha+\frac2q)}+\delta\left(\int_{\mathbb{R}^2}|w|^{2q}\eta^{2\ell}dx\right)\\
%&&+ C(\delta,q,\ell)R^2+\delta\left(\int_{\mathbb{R}^2}|h|^{2q}\eta^{2\ell}dx\right)
\eeno
Thus for the term $I_5$, Young inequality implies that
\beno
I_5&\leq&C(q)R^{2\beta} \int_{\mathbb{R}^2} |\nabla u||\nabla b||h|^{2q-3}(\eta^{2\ell})dx\\
%&\leq& C(q,\ell)R^{\frac1q+2\beta}\left(R^{-1+\alpha+\frac1q}  +\left(\int_{\mathbb{R}^2}|w|^{2q}\eta^{2\ell}dx\right)^{\frac{1}{2q}} \right)\cdot\left(R^{-1+\beta+\frac1q}  +\left(\int_{\mathbb{R}^2}|h|^{2q}\eta^{2\ell}dx\right)^{\frac{1}{2q}} \right) \\
%&&\cdot \left(\int_{\mathbb{R}^2}h^{2q}\eta^{2\ell}dx\right)^{\frac{2q-3}{2q}}\\
&\leq& C(\delta,q,\ell)R^{\frac{2q}{3}(-2+\alpha+3\beta+\frac3q)}+ C(\delta,q,\ell)R^{q(-1+\alpha+2\beta+\frac2q)}+\delta\left(\int_{\mathbb{R}^2}|w|^{2q}\eta^{2\ell}dx\right)\\
&&+ C(\delta,q,\ell)R^{2+4\beta q}+\delta\left(\int_{\mathbb{R}^2}|h|^{2q}\eta^{2\ell}dx\right)
\eeno
where $\delta>0$, to be decided.

{\bf Estimate of $I_1$.} Noting $\ell\geq q$, by Young inequality we have
\beno
I_1
&=& C(\ell,q)R^{2\alpha-2}\int_{\mathbb{R}^2} w^{2q-2}(\eta^{2\ell-2})dx+C(\ell,q)R^{3\alpha-1}\int_{\mathbb{R}^2} w^{2q-2}(\eta^{2\ell-1})dx\\
&&+C(\ell,q)R^{\alpha-1}\int_{\mathbb{R}^2} w^{2q-1}\eta^{2\ell-1}dx=I_{11}+\cdots+I_{13},
\eeno
where
\beno
I_{11}&\leq& \delta \int_{\mathbb{R}^2}w^{2q}\eta^{(2\ell-2)\frac{q}{q-1}}dx+C(\delta,\ell,q)R^{2+q(2\alpha-2)},
\eeno
\beno
I_{12}&\leq& \delta \int_{\mathbb{R}^2}w^{2q}\eta^{(2\ell-1)\frac{q}{q-1}}dx+C(\delta,\ell,q)R^{2+q(3\alpha-1)},
\eeno
and
\beno
I_{13}&\leq& \delta \int_{\mathbb{R}^2}w^{2q}\eta^{(2\ell-1)\frac{2q}{2q-1}}dx+C(\delta,\ell,q)R^{2+2q(\alpha-1)}.
\eeno

Similarly, for the term $I_2$, we get
\beno
I_{2}&\leq& \delta \int_{\mathbb{R}^2}(w^{2q}+h^{2q})\eta^{(2\ell)}dx+C(\delta,\ell,q)R^{2+q(2\alpha+2\beta)}\\
&&+ \delta \int_{\mathbb{R}^2}(w^{2q}+h^{2q})\eta^{(2\ell-1)\frac{q}{q-1}}dx+C(\delta,\ell,q)R^{2+q(2\alpha-1+\beta)}
\eeno

{\bf Estimate of $I_3$.} By H\"{o}lder and Young inequalities we have
\beno
I_{3}&\leq& \delta \int_{\mathbb{R}^2}h^{2q}\eta^{(2\ell-2)\frac{q}{q-1}}dx+C(\delta,\ell,q)R^{2+q(2\beta-2)}\\
&&+ \delta \int_{\mathbb{R}^2}h^{2q}\eta^{(2\ell-1)\frac{q}{q-1}}dx+C(\delta,\ell,q)R^{2+q(\alpha-1+2\beta)}
\\&&+ \delta \int_{\mathbb{R}^2}h^{2q}\eta^{(2\ell-1)\frac{2q}{2q-1}}dx+C(\delta,\ell,q)R^{2+2q(\beta-1)}
\eeno
Similarly,
\beno
I_{4}&\leq& \delta \int_{\mathbb{R}^2}(w^{2q}+h^{2q})\eta^{(2\ell)}dx+C(\delta,\ell,q)R^{2+q(4\beta)}\\
&&+ \delta \int_{\mathbb{R}^2}(w^{2q}+h^{2q})\eta^{(2\ell-1)\frac{q}{q-1}}dx+C(\delta,\ell,q)R^{2+q(-1+3\beta)}
\eeno

Hence, firstly taking $\ell=q$ and
 $\delta<\frac1{32}$; secondly, for fixed $\alpha<\frac13$ with $\beta<-\alpha$, we take the minimum $q_0$ satisfying the following conditions
\beno
2+q(2\alpha-2)\leq  0,\quad 2+q(3\alpha-1)\leq  0,\quad 2+2q(\alpha-1)\leq 0,
\eeno
and
\beno
2+4\beta q\leq 0,\quad 2+q(2\alpha+2\beta)\leq 0.
\eeno
Obviously, $q_0$ is as in (\ref{eq:q0}). And for any $q>q_0$, we %have
%\beno
%2+q(2\alpha-2)< 0,\quad 2+q(3\alpha-1)<  0,\quad 2+2q(\alpha-1)< 0,
%\eeno
%and
%\beno
%2+4\beta q<0,\quad 2+q(2\alpha+2\beta)< 0.
%\eeno
write
\beno
\gamma_0=\max\{2+q(3\alpha-1), 2+4\beta q, 2+q(2\alpha+2\beta) \}<0.
\eeno
Then we get
\beno
\int_{\mathbb{R}^2}w^{2q}\eta^{2\ell}+h^{2q}\eta^{2\ell}dx&\leq& C(\ell,q)\left[R^{2+q(2\alpha-2)}+R^{2+q(3\alpha-1)}+ R^{2+2q(\alpha-1)}\right]\\
&&+C(q,\ell)R^{2+4\beta q}+C(\ell,q)R^{2+q(2\alpha+2\beta)}
\eeno
Choose $R=2^{k+1}$ with $k\in \mathbb{N}$ such that
\beno
\int_{2^k\leq |x|\leq 2^{k+1}}w^{2q}+h^{2q}dx\leq C(\alpha,\beta,q) 2^{k\gamma_0}
\eeno
Consequently, we get
\ben\label{eq:estimate in exterior}
\int_{\mathbb{R}^2\setminus{B_1}}w^{2q}+h^{2q}dx\leq C(\alpha,\beta,q,c_0)<\infty,
\een
for any $q>q_0.$

{\bf Arguments for the estimate in $B_1$.} Firstly,
\beno
\int_{\mathbb{R}^2}w^{2q}+h^{2q}dx<\infty,\quad q>q_0,
\eeno
due to the regularity of the solutions. Secondly, by (\ref{eq:CZ estimate}) we have
\beno
\int_{B_R}|\nabla u|^{2q}+|\nabla b|^{2q}dx\leq C(q)\int_{\mathbb{R}^2}w^{2q}+h^{2q}dx+ C(q)R^{-2q}\int_{B_{2R}}(|u|+|b|)^{2q}dx,
\eeno
and thus
\ben\label{eq:gradient and vorticity}
\int_{\mathbb{R}^2}|\nabla u|^{2q}+|\nabla b|^{2q}dx\leq C(q)\int_{\mathbb{R}^2}w^{2q}+h^{2q}dx<\infty,\quad q>q_0,
\een
where we used the growth estimates (\ref{eq:growth at infty-u}) and  (\ref{eq:growth at infty-b}). Finally, for
the cut-off function $\eta_1$, i.e.
\begin{align*} \eta_1(x)=\left\{
\begin{aligned}
&1,\quad x\in B_1,\\
&0, \quad x\in B_{2}^c,
\end{aligned}
\right. \end{align*}
one can also obtain that
(\ref{eq:w estimate-mhd-general}) and
\beno
\int_{\mathbb{R}^2}w^{2q}\eta_1^{2\ell}+h^{2q}\eta_1^{2\ell}dx&\leq&C(q)\int_{\mathbb{R}^2}(|\nabla u|+|\nabla b|)^{2q-2}\eta_1^{2\ell-2}+(|\nabla u|+|\nabla b|)^{2q-1}\eta_1^{2\ell-1}dx,
\eeno
which can be controlled by
\beno
C(q)\int_{\mathbb{R}^2}(|\nabla (u\eta_1)|+|\nabla (b\eta_1)|)^{2q-2}+(|\nabla (u\eta_1)|+|\nabla (b\eta_1)|)^{2q-1}dx+C(q).\eeno
Using (\ref{eq:CZ estimate}) again, we get
\beno
\int_{\mathbb{R}^2}w^{2q}\eta_1^{2\ell}+h^{2q}\eta_1^{2\ell}dx\leq  \frac12\int_{\mathbb{R}^2}w^{2q}\eta_1^{2\ell}+h^{2q}\eta_1^{2\ell}dx+C(q),
\eeno
which and (\ref{eq:estimate in exterior}) imply that
\beno
\int_{\mathbb{R}^2}w^{2q}+h^{2q}dx\leq C(\alpha,\beta,q,c_0)<\infty,
\eeno
for any $q>q_0.$
And the required inequality follows by using (\ref{eq:CZ estimate}), (\ref{eq:growth at infty-u}) and  (\ref{eq:growth at infty-b}) again.

Thus the proof of Proposition \ref{thm:generalization of KNSS2} is complete.

\begin{lem}\label{thm:generalization of KNSS3}
Let $(u,b,\pi)$ be a smooth solution of the 2D MHD equations (\ref{eq:MHD-2D}) defined over the entire plane satisfying the growth estimates (\ref{eq:growth at infty-u}) with $0<\alpha<\frac13$.
Moreover, we assume that $b$ satisfies (\ref{eq:growth at infty-b}) with $\beta<-\alpha$.
Then
\beno
\|\nabla (|w|^{q-1})\|_{L^q(\mathbb{R}^2)}+\|\nabla (|h|^{q-1})\|_{L^2(\mathbb{R}^2)}\leq C(\alpha,\beta,q,c_0)<\infty,
\eeno
where
\beno
q>q_0+1=\max\{\frac{2}{1-3\alpha}, \frac{-1}{\alpha+\beta},-\frac{1}{2\beta}\}+1.
\eeno
\end{lem}

{\bf Proof of Lemma \ref{thm:generalization of KNSS3}.}
On the other hand, let $\phi(x)\in C_0^\infty(B_R)$ and $0\leq \phi\leq 1$ satisfying
\begin{align*} \phi(x)=\left\{
\begin{aligned}
&1,\quad x\in B_R,\\
&0, \quad x\in B_{2R}^c
\end{aligned}
\right. \end{align*}

Using similar estimates as in (\ref{eq:energy estimate-w}) and (\ref{eq:energy estimate-h}),
multiply $\phi^{2q}w^{2q-4}w$ and $\phi^{2q}h^{2q-4}h$ on both sides of (\ref{eq:MHDv-2D}) with $q>2$, and we have
\ben\label{eq:energy estimate-w'}
I'&\doteq&(2q-4)\int_{\mathbb{R}^2}|\nabla w|^2 w^{2q-4}\phi^{2q}dx\nonumber\\
&\leq&\frac{1}{2q-2}\int_{\mathbb{R}^2} w^{2q-2}\triangle (\phi^{2q})dx+\frac{1}{2q-2}\int_{\mathbb{R}^2} w^{2q-2}u\cdot \nabla(\phi^{2q})dx\nonumber\\
&&+C(q)\int_{\mathbb{R}^2}h^2~ |b|^2w^{2q-4}(\phi^{2q})dx -\int_{\mathbb{R}^2}w^{2q-4}w h~ b\cdot \nabla (\phi^{2q})dx\nonumber\\
&\doteq&I_1'+\cdots+I_4',
\een
and
\ben\label{eq:energy estimate-h'}
II'&\doteq&(2q-4)\int_{\mathbb{R}^2}|\nabla h|^2 h^{2q-4}\phi^{2q}dx\nonumber\\
&\leq&\frac{1}{2q-2}\int_{\mathbb{R}^2} h^{2q-2}\triangle (\phi^{2q})dx+\frac{1}{2q-2}\int_{\mathbb{R}^2} h^{2q-2}u\cdot \nabla(\phi^{2q})dx\nonumber\\
&&-\int_{\mathbb{R}^2}  w ~h^{2q-4}hb\cdot \nabla(\phi^{2q})dx+C(q)\int_{\mathbb{R}^2}w^2~ |b|^2h^{2q-4}(\phi^{2q})dx \nonumber\\
&&+C\int_{\mathbb{R}^2} |\nabla u||\nabla b||h|^{2q-3}(\phi^{2q})dx\nonumber\\
&\doteq&II_1'+\cdots+II_5'
\een

Since
\beno
\|w\|_{2q} +\|h\|_{2q}   <\infty,
\eeno
for any $q>q_0$,
we have
\beno
 \|\nabla(|w|^{q-1})\|^2_{L^2(B_R)} + \|\nabla(|h|^{q-1})\|_{L^2(B_R)}^2 \leq C\int_{B_{2R}}|w|^{2q-2}+|h|^{2q-2}+|\nabla u||h|^{2q-2} dx<\infty,
\eeno
for any $q>q_0+1$.
Then the proof is complete. \endproof

{\bf Proof of  Theorem \ref{thm:generalization of KNSS}:}
For $q>q_0+1$, we still consider the inequalities (\ref{eq:energy estimate-w'}) and (\ref{eq:energy estimate-h'}).
Now we estimate the term $I_3'$, since
\beno
I'_{3}\leq C(q) \int_{\mathbb{R}^2}|b|^2w^{2q-2}(\phi^{2q})dx+\int_{\mathbb{R}^2}|b|^2h^{2q-2}(\phi^{2q})dx\doteq I'_{31}+I'_{32},
\eeno
where
\beno
I'_{31}&=&\int_{\mathbb{R}^2}|b|^2w^{2q-2}(\phi^{2q})dx\\
&\leq& \left(\int_{\mathbb{R}^2}|b|^{2p'}dx\right)^{\frac{1}{p'}} \left(\int_{\mathbb{R}^2}|\tilde{w}|^{2p}dx\right)^{\frac{1}{p}}\\
&\leq& C(q)\|b\|_{L^{p'}(\mathbb{R}^2)} \|b\|_\infty\|\tilde{w}\|_{\frac{2q}{q-1}}^\theta\|\nabla\tilde{w}\|_{2}^{2-\theta}
\eeno
where $\tilde{w}=|w|^{q-1}$ and we used H\"{o}lder inequality, Lemma \ref{thm:generalization of KNSS3}, and Gagliardo-Nirenberg inequality( for example, see Lemma II.3.3 in \cite{Galdi}). Let $p=4q$, then
\beno
\theta=\frac{2q}{p(q-1)}=\frac{1}{2(q-1)}
\eeno
Taking $p_1=8q-2$, noting that (\ref{eq:gradient and vorticity}), by Gagliardo-Nirenberg inequality we have
\beno
I'_{31}&\leq& C(q)\|b\|_{L^{p'}(\mathbb{R}^2)} \|b\|_{L^1(\mathbb{R}^2)}^{(\frac{p_1-2}{3p_1-2})}\|h\|_{L^{p_1}(\mathbb{R}^2)}^{(\frac{2p_1}{3p_1-2})}\|\tilde{w}\|_{\frac{2q}{q-1}}^\theta\|\nabla\tilde{w}\|_{2}^{2-\theta}\\
&\leq&C(q) \|b\|_{L^{p'}(\mathbb{R}^2)} \|b\|_{L^1(\mathbb{R}^2)}^{(\frac{p_1-2}{3p_1-2})}\|\tilde{h}\|_{L^{\frac{2q}{q-1} }(\mathbb{R}^2)}^{\frac{4q}{q-1} \frac{1}{3p_1-2}}\|\nabla\tilde{h}\|_{2}^{ (\frac{p_1-2q}{p_1(q-1)})(\frac{2p_1}{3p_1-2})}
\|\tilde{w}\|_{\frac{2q}{q-1}}^\theta\|\nabla\tilde{w}\|_{2}^{2-\theta}\\
&\leq& C(q)\|b\|_{L^{p'}(\mathbb{R}^2)} \|b\|_{L^1(\mathbb{R}^2)}^{(\frac{p_1-2}{3p_1-2})}\|\tilde{h}\|_{L^{\frac{2q}{q-1} }(\mathbb{R}^2)}^{\frac{4q}{q-1} \frac{1}{3p_1-2}}\|\nabla\tilde{h}\|_{2}^{ \theta}
\|\tilde{w}\|_{\frac{2q}{q-1}}^\theta\|\nabla\tilde{w}\|_{2}^{2-\theta},
\eeno
since
\ben\label{eq:theta}
(\frac{p_1-2q}{p_1(q-1)})(\frac{2p_1}{3p_1-2})=\theta.
\een
Due to (\ref{eq:growth at infty-b}) and Proposition \ref{thm:generalization of KNSS2}, we have
\beno
\|b\|_{L^{p'}(\mathbb{R}^2)}\leq C(\beta)\|b\|_{L^{1}(\mathbb{R}^2)}^{\frac{1}{p'}}
\eeno
hence there exists a positive number $\|b\|_{L^{1}(\mathbb{R}^2)}=\varepsilon_2=\varepsilon_2(\alpha,\beta,q,c_0)$ such that
\beno
I'_{31}
&\leq& C(q,\alpha,\beta,c_0) \|b\|_{L^1(\mathbb{R}^2)}^{(\frac{p_1-2}{3p_1-2})+\frac{1}{p'}}\|\nabla\tilde{h}\|_{2}^{ \theta}
\|\nabla\tilde{w}\|_{2}^{2-\theta}\\
&\leq& \frac{1}{16}[\|\nabla\tilde{h}\|_{2}^2+ \|\nabla\tilde{w}\|_{2}^2]
\eeno

The term $I'_{32}$ and $II_4'$ are similar, hence we have
\beno
I'_{3}+II_4'&=&\int_{\mathbb{R}^2}w^2~ |b|^2h^{2q-2}(\phi^{2q})dx\leq \frac{1}{8}[\|\nabla\tilde{h}\|_{2}^2+ \|\nabla\tilde{w}\|_{2}^2]
\eeno

Next we estimate the term $II_5'$.  Noting that (\ref{eq:gradient and vorticity}), using Gagliardo-Nirenberg inequality again we have
%\beno
%II_5'=\int_{\mathbb{R}^2} |\nabla u||\nabla b||h|^{2q-1}(\phi^{2q})dx
%\eeno
\beno
II_5'&=&\int_{\mathbb{R}^2}|\nabla b||\nabla u|h^{2q-3}(\phi^{2q})dx\\
&\leq &\|\nabla u\|_{2p(q-1)}\|\nabla b\|_{\frac{2p(q-1)(2q-2)}{2pq-2p-1}}^{2q-2}\\
&\leq& C(q)\left(\int_{\mathbb{R}^2}|h|^{p'}dx\right)^{\frac{1}{p'}} \left(\int_{\mathbb{R}^2}|h^{2q-2}|^{p}dx\right)^{\frac{2q-3}{p(2q-2)}}\|\nabla u\|_{2p(q-1)}\\
&\leq& C(q)\|h\|_{L^{p'}(\mathbb{R}^2)} \left[\|\tilde{w}\|_{\frac{2q}{q-1}}^\theta\|\nabla\tilde{w}\|_{2}^{2-\theta} +\|\tilde{h}\|_{\frac{2q}{q-1}}^\theta\|\nabla\tilde{h}\|_{2}^{2-\theta}\right]
\eeno
where $p=4q$ and
\beno
\theta=\frac{2q}{p(q-1)}=\frac{1}{2(q-1)}
\eeno
Taking $p_1=8q-2$ and $\gamma=\frac{p_1-2}{3p_1-2}$ we have
\beno
C(q)\|h\|_{L^{p'}(\mathbb{R}^2)}&\leq& C(q)\left(\int_{\mathbb{R}^2} |h|^{\frac{p'\gamma p_1}{p_1-p'(1-\gamma)}}dx\right)^{\frac{p_1-p'(1-\gamma)}{p_1p'}} \|h\|_{L^{p_1}(\mathbb{R}^2)}^{(\frac{2p_1}{3p_1-2})}\\
&\leq&C(q) \|h\|_{L^{\frac{p'\gamma p_1}{p_1-p'(1-\gamma)}}(\mathbb{R}^2)}^{\gamma} \|\tilde{h}\|_{L^{\frac{2q}{q-1} }(\mathbb{R}^2)}^{\frac{4q}{q-1} \frac{1}{3p_1-2}}\|\nabla\tilde{h}\|_{2}^{ (\frac{p_1-2q}{p_1(q-1)})(\frac{2p_1}{3p_1-2})},
\eeno
and since
\beno
\frac{p'\gamma p_1}{p_1-p'(1-\gamma)}=\frac{p_1\frac{p}{p-1}\cdot \frac{p_1-2}{3p_1-2} }{p_1-\frac{p}{p-1}\cdot\frac{2p_1}{3p_1-2}}=\frac{4q(8q-4)}{(24q-8)(4q-1)-8q}=\frac{4q^2-2q}{12q^2-8q+1}
\eeno
thence by (\ref{eq:theta}) we have
\beno
II_5'
&\leq& C(q)\||h|^{\frac13}\|_{L^{1}(\mathbb{R}^2)}^{3\gamma\frac{6q-2}{6q-1}}\|h\|_{L^{2q}(\mathbb{R}^2)}^{\frac{\gamma}{6q-1}}\|\tilde{h}\|_{L^{\frac{2q}{q-1} }(\mathbb{R}^2)}^{\frac{4q}{q-1} \frac{1}{3p_1-2}}\|\nabla\tilde{h}\|_{2}^{ \theta}\\
&&\cdot \left[\|\tilde{w}\|_{\frac{2q}{q-1}}^\theta\|\nabla\tilde{w}\|_{2}^{2-\theta} +\|\tilde{h}\|_{\frac{2q}{q-1}}^\theta\|\nabla\tilde{h}\|_{2}^{2-\theta}\right],
\eeno
where we used H\"{o}lder inequality, since
\beno
\frac13<\frac{4q^2-2q}{12q^2-8q+1}<2q
\eeno
for $q>2$. Hence there exists a positive number $\||h|^{\frac13}\|_{L^{1}(\mathbb{R}^2)}\leq \varepsilon_3(\alpha,\beta,q, c_0)$ such that
\beno
II_5'
&\leq& \frac{1}{16}[\|\nabla\tilde{h}\|_{2}^2+ \|\nabla\tilde{w}\|_{2}^2]
\eeno

Recalling the inequalities (\ref{eq:energy estimate-w'}) and (\ref{eq:energy estimate-h'}), using
 the growth (\ref{eq:growth at infty-u}), (\ref{eq:growth at infty-b}) and the above estimates, by  Proposition \ref{thm:generalization of KNSS2} and  Lemma \ref{thm:generalization of KNSS3} we get
\beno
&&\int_{\mathbb{R}^2}|\nabla w|^2 w^{2q-4}\phi^{2q}dx+\int_{\mathbb{R}^2}|\nabla h|^2 h^{2q-4}\phi^{2q}dx\\
&\leq &C(\alpha,\beta,q,c_0)[R^{-2}+ R^{\alpha-1}+  R^{\beta-1}]
\eeno
and $R\rightarrow\infty$ implies that
\beno
\nabla(|w|^{q-1})\equiv 0,\quad \nabla(|h|^{q-1})\equiv 0,
\eeno
which yields that
\beno
w\equiv C,\quad h\equiv C,
\eeno
and it follows from Proposition \ref{thm:generalization of KNSS2} that $C\equiv 0$ and $u,b$ are constants.

The proof of  Theorem \ref{thm:generalization of KNSS} is complete by taking $\varepsilon_0=\min\{\varepsilon_2,\varepsilon_3\}$.

%{\bf Step I. $w,h\in L^q(\mathbb{R}^2)$} for some large $q$, hence
%\beno
%\nabla u,\nabla b\in L^q(\mathbb{R}^2)
%\eeno
%
%{\bf Step II: the key point!???}
%\beno
% u\in L^\infty(\mathbb{R}^2)
%\eeno
%
%{\bf Step III. $u,b\in L^\infty(\mathbb{R}^2)$} and
%\beno
%\|b\|_{L^1(\mathbb{R}^2)}\|b\|_{L^\infty(\mathbb{R}^2)}\leq C_*\min\{\mu\nu, \mu^{\frac12}\nu^{\frac32}\}
%\eeno
%
%
%The difficulty: $\alpha<\frac17$, no maximum principle; here $\nabla u\in L^q$(large q) is not easy, $q\leq 6$ is OK;
%
%
%
%!!!:Let $\alpha\leq \frac29$, then we have $\nabla u\in L^6$, which implies $u=b=C$ by W-W theorem 1.1.

\section{Proof of  Theorem \ref{thm:generalization of KNSS'}}

\begin{prop}\label{thm:generalization of KNSS''}
Let $(u,b,\pi)$ be a smooth solution of the 2D MHD equations (\ref{eq:MHD-2D}) defined over the entire plane satisfying the growth estimates $\nabla u\in L^{q_0}(\mathbb{R}^2)$ for $2<q_0<\infty$, $\nabla b\in L^{\infty}(\mathbb{R}^2)$ and
\beno
\|b\|_{L^1(\mathbb{R}^2)}+\||h|^{1/3}\|_{L^1(\mathbb{R}^2)}\leq \varepsilon_1,
\eeno
where $\varepsilon_1$ is sufficiently small depending on $q_0,$ $\|\nabla u\|_{L^{q_0}}$ and $\|\nabla b\|_{L^\infty}$. Then
\beno
\nabla u\in L^{p}(\mathbb{R}^2), \nabla b\in L^{p}(\mathbb{R}^2),
\eeno
for any $p\geq q_0.$
\end{prop}

{\bf Proof of Proposition \ref{thm:generalization of KNSS''}.} By Lemma \ref{lem:Galdi}, there exists $R>0$ such that
\ben\label{eq:bound of nabla u}
|u(x)|\leq (1+|x|)^{\frac{q_0-2}{q_0}},\quad |x|>R,
\een
since $\nabla u\in L^{q_0}(\mathbb{R}^2)$, and we also have $b(x)\in L^p(\mathbb{R}^2)$ with $1\leq p\leq \infty$ by Gagliardo-Nirenberg inequality satisfying
\beno
\|b\|_{L^p(\mathbb{R}^2)}\leq C(p, \|\nabla b\|_{\infty})
\eeno
%\beno
%|b(x)|\leq C(1+|x|)^{\frac{\frac{q_0^2}{q_0-1}-2}{\frac{q_0^2}{q_0-1}}}
%\eeno
Moreover, using (\ref{eq:gradient and vorticity}) again, we have $\nabla b(x)\in L^p(\mathbb{R}^2)$ with $1<p\leq \infty$ and
\beno
\|\nabla b\|_{L^p(\mathbb{R}^2)}\leq C(p, \|\nabla b\|_{\infty}).
\eeno

Recalling the inequalities (\ref{eq:energy estimate-w'}) and (\ref{eq:energy estimate-h'}) with $q-1=\frac{q_0}{2}$, we have
\ben\label{eq:w q0}
 &&\|\nabla(|w|^{q-1})\|^2_{L^2(B_R)} + \|\nabla(|h|^{q-1})\|_{L^2(B_R)}^2 \nonumber\\
 &\leq& C(q_0,\|\nabla b\|_{\infty})\int_{B_{2R}\setminus{B_R}}(R^{-\frac{2}{q_0}}+|b|R^{-1})(|w|^{q_0}+|h|^{q_0})dx\nonumber\\
  &&+C(q_0,\|\nabla b\|_{\infty})\int_{B_{2R}}|b|^2(|w|^{q_0}+|h|^{q_0})+|\nabla u||\nabla b|^{q_0} dx<\infty,
\een
Thus
\beno
\nabla(|w|^{\frac{q_0}{2}}),\nabla(|h|^{\frac{q_0}{2}})\in L^2(\mathbb{R}^2),
\eeno
which implies
Proposition \ref{thm:generalization of KNSS''}.

{\bf Proof of  Theorem \ref{thm:generalization of KNSS'}:} One can make the same argyuments with the two terms of $I_{31}'$ and $II_5'$ as in the proof of Theorem \ref{thm:generalization of KNSS} by noting that (\ref{eq:w q0}) and (\ref{eq:bound of nabla u}).

For $1<q_0\leq 2$, we refer to Theorem 1.1 and 1.2 of \cite{WW}, which is an immediate corollary.
The proof is complete.

\noindent {\bf Acknowledgments.}
W. Wang was supported by NSFC under grant 11671067 and
 "the Fundamental Research Funds for the Central Universities".

\end{document}